# FINITELY PRESENTED SUBGROUPS OF AUTOMATIC GROUPS AND THEIR ISOPERIMETRIC FUNCTIONS

GILBERT BAUMSLAG, MARTIN R. BRIDSON,
CHARLES F. MILLER III AND HAMISH SHORT

ABSTRACT. We describe a general technique for embedding certain amalgamated products into direct products. This technique provides us with a way of constructing a host of finitely presented subgroups of automatic groups which are not even asynchronously automatic. We can also arrange that such subgroups satisfy, at best, an exponential isoperimetric inequality.

**1. Introduction.** Despite the intensive study of automatic and asynchronously automatic groups over the past few years, very little is known about their finitely presented subgroups. In fact, until now, the only known finitely presented subgroup of an automatic group which is not automatic has been an example of John Stallings [18]. This example was pointed out by G. Mess (cited in [9]), who noted that Stallings' group $S$ is a subgroup of a direct product $D$ of three free groups of rank two – this observation actually goes back to R. Bieri (see remark 4.3 of [5]). The point here is that finitely generated free groups are automatic and the direct product of finitely many automatic groups is again automatic, so $D$ is an automatic group. Since automatic groups are of type $FP_\infty$ and $S$ is not even of type $FP_3$, it follows that $S$ is not automatic. (We refer the reader to [8] for definitions and standard facts concerning automatic groups.)

Our basic method for constructing examples of finitely presented subgroups of automatic groups which are not automatic was motivated by this example of Stallings. It centres on the notion of doubling a group along a subgroup. In order to explain, let $A$ be a group and let $\overline{A}$ be an isomorphic copy of $A$, where the underlying isomorphism $\theta$ maps $a \in A$ to $\overline{a} \in \overline{A}$ ; if $C$ is a subgroup of $A$ then we denote its image in $\overline{A}$ by $\overline{C}$. We then term an amalgamated product
$$G = \{A * \overline{A}; C = \overline{C}\},$$
where $C$ is identified with $\overline{C}$ according to the map $\theta$, the *double of $A$ along $C$.*

*Key words and phrases.* Automatic and asynchronously automatic groups, isoperimetric functions, finitely presented groups, subgroups.

This work was supported in part by NSF Grants #9103098 (first author) and #9203500 (second author).

There is an analogous notion for HNN extensions. If $A$ is a group and $C$ is a subgroup of $A$ we write
$$A*_C = \{A, s : s^{-1}cs = c, \forall c \in C\}$$
which is the HNN extension obtained by adding a new stable letter $s$ which commutes with the subgroup $C$. As is well known, the subgroup of $A*_C$ generated by $A$ and $s^{-1}As$ is their free product with amalgamation over $C$. Hence $A*_C$ contains an isomorphic copy of the double of $A$ along $C$.

Our constructions rely on the following basic embedding result.

**Theorem 1.** *Let $A$ be a group with normal subgroup $C$ and let $G$ be the double of $A$ along the normal subgroup $C$. Then $G$ embeds in*
$$D = (A/C * \overline{A/C}) \times A$$
*the direct product of $A$ with the ordinary free product of two copies of $A/C$. Similarly, the HNN extension $A*_C$ of $A$ by a stable letter $s$ commuting with $C$ embeds in*
$$D' = (A/C * \langle s \rangle) \times A$$
*the direct product of $A$ with the ordinary free product of $A/C$ with the infinite cyclic group $\langle s \rangle$ generated by $s$.*

As an indication of how we make use of these embeddings, recall that if $H$ and $H'$ are automatic (respectively biautomatic) groups then the groups $H * H'$ and $H \times H'$ are also automatic (respectively biautomatic). It follows from the above embedding result that, if $C$ is normal in $A$ and both $A$ and $A/C$ are (bi)automatic, then both $A*_C$ and the double $G$ of $A$ along $C$ can be embedded in (bi)automatic groups.

Three kinds of examples will be described in detail below. All of our examples are obtained from this and similar embedding theorems by choosing the group $A$ and its subgroup $C$ suitably. We can then, for instance, make use of the following.

**Theorem 2.** *Let $1 \to C \to A \to Q \to 1$ be a short exact sequence of finitely generated infinite groups. Suppose that $A$ is word hyperbolic and $Q$ is biautomatic. Then the double $G$ of $A$ along $C$ satisfies, at best, an exponential isoperimetric inequality.*

By combining Theorem 2 with our embedding results, we construct examples of the following type.

**Theorem 3.** *There exists a biautomatic group $B$ and a finitely presented subgroup $G \leq B$ such that $G$ is not of type $FP_3$ and its isoperimetric function is strictly exponential. Moreover one can arrange for $B$ to be the fundamental group of a closed manifold of non-positive curvature.*

Since asynchronously automatic groups are of type $FP_3$ [10] the subgroups mentioned in this theorem are not even asynchronously automatic. Stallings' group provides another such example; however it satisfies a cubic isoperimetric inequality (as we shall show in section §6 below).

Our examples serve to show that the difficulty in solving the word problem in a group $G$ is not well measured by its isoperimetric function alone. Indeed as noted above we have, in particular, an example of a finitely presented subgroup $G$ of a biautomatic group $B$ which satisfies a strictly exponential isoperimetric inequality. The word problem in an automatic group can be solved in quadratic time [8]. Thus we can decide whether any given word $w$ in the generators of $G$ is equal to the identity in quadratic time by simply re-expressing $w$ as a word in the generators of $B$.

The techniques which we develop in this article can be applied to the study of isoperimetric inequalities for finitely presented subgroups of finitely presented groups in general. Previous descriptions of groups where the isoperimetric function of a finitely presented subgroup is appreciably larger than that of the ambient group have involved delicate analytic arguments (see [16]). Thus our examples appear to be of interest from a purely geometric perspective. After learning about our work, S.M. Gersten has further studied the isoperimetric functions for subgroups and what he terms "area distortion" [11], and has produced several interesting examples.

It is also worth noting that it is not known whether finitely presented subgroups of word hyperbolic groups can have an isoperimetric function which is not linearly bounded, so it may even be that finitely presented subgroups of word hyperbolic groups are themselves word hyperbolic. (Recall that a finitely presented group is word hyperbolic if and only if it satisfies a linear isoperimetric inequality [16].) Bestvina has recently shown that Gromov's proposed counter–example fails [3], and S.M. Gersten has established that the statement is true for hyperbolic groups of cohomological dimension 2, and for those with diagrammatically aspherical finite presentations [12].

**2. Some elementary embedding facts.** We begin by recording a general embedding observation whose proof is immediate.

**Lemma 4.** *Let $\phi : H \to K$ and $\psi : H \to L$ be group homomorphisms. Then the homomorphism from $H$ to $K \times L$ defined by $h \mapsto (\phi(h), \psi(h))$ is an embedding if and only if $\ker \phi \cap \ker \psi = \{1\}$.*

Theorem 1 is a special case of the following two easy results which provide the basic embedding facts we need.

**Theorem 5.** *Let $A$ be a group with subgroups $N \leq C$ where $N$ is normal in $A$. Let $G$ be the double of $A$ along the subgroup $C$, that is $G = \{A * \overline{A}; C = \overline{C}\}$. Then $G$ embeds in the group*

$$D = \{A/N * \overline{A}/\overline{N}; C/N = \overline{C}/\overline{N}\} \times A$$

*which is the direct product of $A$ with the double of $A/N$ along $C/N$. In particular, if $N = C$ so that $C$ is normal, $G$ embeds in*

$$D = (A/C * \overline{A}/\overline{C}) \times A$$

*the direct product of $A$ with the ordinary free product of two copies of $A/C$.*

*Proof.* Let $\phi : G \to \{A/N * \overline{A}/\overline{N} : C/N = \overline{C}/\overline{N}\}$ be the quotient mapping whose kernel is the normal subgroup $N$. Let $\psi : G \to A$ be the homomorphism defined by $\psi(a) = a$ for all $a \in A$ and $\psi(\overline{a}) = a$ for all $\overline{a} \in \overline{A}$. Then $\ker \psi \cap A = \{1\}$ and so $\ker \phi \cap \ker \psi = \{1\}$. The result now follows from Lemma 4. □

Here is the analogous result for HNN extensions:

**Theorem 6.** *Let $A$ be a group with subgroups $N \subseteq C$ where $N$ is normal in $A$. Then the HNN extension $A*_C$ of $A$ by a stable letter $s$ commuting with $C$ embeds in the group*

$$D' = (A/N*_{C/N}) \times A$$

*which is the direct product of $A$ with $A/N*_{C/N}$. In particular, if $N = C$ so that $C$ is normal, $G$ embeds in*

$$D' = (A/C * \langle s \rangle) \times A$$

*the direct product of $A$ with the ordinary free product of $A/C$ with the infinite cyclic group $\langle s \rangle$ generated by $s$.*

*Proof.* Let $\phi : A*_C \to A/N*_{C/N}$ be the quotient map whose kernel is the normal subgroup $N$. Let $\psi : A*_C \to A$ be the map which is the identity on $A$ and sends the commuting stable letter to $1 \in A$. Then $\ker \psi \cap A = \{1\}$ and so $\ker \phi \cap \ker \psi = \{1\}$. The result now follows from Lemma 4. □

We define $d(A)$ to be the minimum number of generators of $A$, $d_A(C)$ to be the minimum number of elements of $C$ which, together with all of their conjugates in $A$, generate $C$ and

$$\omega_A(C) = \begin{cases} 0, & \text{if } C \text{ is finitely generated} \\ d(A/C) & \text{if } C \text{ is not finitely generated.} \end{cases}$$

**Theorem 7.** *Let $C$ be a normal subgroup of the finitely presented group $A$ and assume $d_A(C) < \infty$. Let $G$ be the double of $A$ along $C$. Suppose that $\{a_1, \ldots, a_m\}$ is a finite set of elements of $A$ whose images in $A/C$ generate $A/C$. Let $F_m$ be the free group with basis $t_1, \ldots, t_m$ and let each $t_i$ act on $A$ by conjugation by the corresponding $a_i$. Similarly, let $t_i$ act on the copy $\overline{A}$ of $A$ by conjugation by the corresponding $\overline{a}_i$. Then the split extension $E = G \rtimes F_m$ is finitely presented and embeds in the finitely presented group*

$$D = ((A/C * \overline{A}/\overline{C}) \rtimes F_m) \times A$$

*the direct product of $A$ with the split extension of the ordinary free product of two copies of $A/C$.*

*Proof.* Since $d_A(C) < \infty$ and the images of $a_1, \ldots, a_m$ generate $A$ there is a finite subset $\{c_1, \ldots, c_k\}$ of $C$ so that $C$ is generated by the conjugates $c_j^h$ of $c_1, \ldots, c_k$ by the elements $h \in gp(a_1, \ldots, a_m)$.

Let $H = A \rtimes F_m$ be the semidirect product of $A$ with the free group $F_m$ of rank $m$ on $t_1, \ldots, t_m$, where the action of $t_i$ on $A$ is conjugation by $a_i$. Clearly $H$ is finitely presented since $A$ is finitely presented and $m$ is finite. Let $B$ be the subgroup of $H$ generated by $C$ and $F_m$. Since $C$ is generated by the conjugates $c_j^h$ of $c_1, \ldots, c_k$ by the elements $h \in gp(a_1, \ldots, a_m)$ and the $t_i$ act as conjugation by $a_i$, it follows that $B$ is finitely generated by $c_1, \ldots, c_k, t_1, \ldots, t_m$.

Let $E_1 = \{H * \overline{H}; B = \overline{B}\}$ be the double of $H$ along the subgroup $B$. Then $E$ is finitely presented since $H$ is finitely presented and $B$ is finitely generated.

Now by examining the natural presentation of $E$ given by the above construction and using the relations $t_i = \overline{t}_i$ to eliminate the generators $\overline{t}_i$, it is clear that

$$E_1 \cong G \rtimes F_m = \{A * \overline{A}; C = \overline{C}\} \rtimes F_m = E$$

where the action of the $t_i$ on $A$ is conjugation by $a_i$ and the action of $t_i$ on $\overline{A}$ is conjugation by $\overline{a_i}$.

If $\phi : E \to E/C$ is the quotient map and $\psi : E \to A$ is the map defined by $\phi(a) = \phi(\overline{A}) = a$ for all $a \in A$ and $\phi(t_j) = a_j$ for $i = 1, \ldots, m$, then as before $\ker \phi \cap \ker \psi = \{1\}$ so the product map embeds $E$ into $D = (E/C) \times A$ by Lemma 4. This completes the proof of Theorem 7. □

Here is the corresponding result for HNN extensions. We omit the proof as it is quite similar to that of the previous theorem.

**Theorem 8.** *Let $C$ be a normal subgroup of the finitely presented group $A$ and assume $d_A(C) < \infty$. Suppose that $\{a_1, \ldots, a_m\}$ is a finite set of elements of $A$ whose images in $A/C$ generate $A/C$. Let $F_m$ be the free group with basis $t_1, \ldots, t_m$ and let each $t_i$ act on $A$ by conjugation by the corresponding $a_i$. Extend this to an action on $A*_C$ by letting each $t_i$ commute with the stable letter $s$. Then the split extension $E' = (A*_C) \rtimes F_m$ of $A*_C$ by $F_m$ is finitely presented and embeds in the finitely presented group*

$$D' = ((A/C * \langle s \rangle) \rtimes F_m) \times A$$

*where $\langle s \rangle$ is the infinite cycle generated by the stable letter $s$ of $A*_C$.*

The following observation is frequently useful.

**Lemma 9.** *Let $\phi : W \to U$ be a homomorphism of groups and let each $w \in W$ act on $U$ as conjugation by $\phi(w)$. Then, the semidirect product $U \rtimes W$ is isomorphic to the direct product $U \times W$.*

*Proof.* In $U \times W$ the elements $(\phi(w), w)$ form an isomorphic copy of $W$ which acts on $U$ by conjugation in the same way $W$ does. Hence $U \times W \cong U \rtimes (\phi(W), W) \cong U \rtimes W$ as required. □

**3. Homological obstructions: examples of Rips' type.** We now turn our attention to the construction of specific classes of examples of finitely presented subgroups of automatic groups that are not themselves automatic. We shall exhibit two types of obstructions in order to prove that the subgroups under consideration are not automatic: homology groups and isoperimetric functions. We begin with homology.

Our first examples make use of a construction due to Rips [17] which some of us considered in [2]. In order to explain, let $\langle X \mid R \rangle$ be a finite presentation of a group $Q$. Then there is a word hyperbolic (indeed small cancellation) group whose presentation takes the form $A = \langle X, a, b \mid S \rangle$, such that $C = gp(a, b)$ is normal in $A$ and such that $A/C \cong Q$. The set $S$ of defining relators is partitioned into two subsets. The first subset consists of relators of the form $r^{-1} w_r(a, b)$ where here $r$ ranges over $R$ and $w_r(a, b)$ denotes a word in $a$ and $b$. The second subset ensures that $C$ is normal in $A$ and consists of relators of the form

$$x_i^{-1} a x_i (u_i(a,b))^{-1}, x_i^{-1} b x_i (v_i(a,b))^{-1}, x_i a x_i^{-1} (w_i(a,b))^{-1}, x_i b x_i^{-1} (z_i(a,b))^{-1},$$

where $x_i$ ranges over $X$. We will refer to such a group $A$ as a *Rips group with (Rips) quotient $Q$ and kernel $C$*. We can arrange that $A$ satisfy as stringent a small cancellation condition as we choose and so, in particular, that it be of cohomological dimension two.

R. Bieri ([5] Theorem B) has proved that a finitely presented normal subgroup of a group of cohomological dimension at most two is either free or of finite index. It follows that $C$ is not finitely related whenever $Q$ is infinite. Indeed it seems likely, although we have not been able to give a complete proof, that $H_2(C)$ is not even finitely generated whenever $Q$ is infinite. These remarks are a prelude to the following observation.

**Lemma 10.** *Let $A$ be a Rips group with kernel $C$ and quotient a biautomatic group $Q$. Assume that $H_2(C)$ is not finitely generated. Then the double $G$ of $A$ along $C$ is a finitely presented subgroup of a biautomatic group, but $G$ is not asynchronously automatic. Similarly, the HNN-extension $A*_C$ is a finitely presented subgroup of a biautomatic group which is not asynchronously automatic.*

*Proof.* It follows immediately from Theorem 1 that $G$ and $A*_C$ are subgroups of a biautomatic groups. Now automatic groups are of type $FP_\infty$ ([8]) while asynchronously automatic groups are of type $FP_3$ ([10]). So it suffices to prove that $G$ and $A*_C$ are not of type $FP_3$. For $G$ this follows immediately from the Mayer-Vietoris sequence for an amalgamated product:

$$\cdots \longrightarrow H_3(A) \oplus H_3(A) \longrightarrow H_3(G) \longrightarrow H_2(C) \longrightarrow H_2(A) \oplus H_2(A) \longrightarrow H_2(G)$$
$$\longrightarrow H_1(C) \longrightarrow H_1(A) \oplus H_1(A) \longrightarrow H_1(G) \longrightarrow \mathbb{Z} \longrightarrow 0.$$

Since $H_3(A) = 0$ and $H_2(A)$ is finitely generated, it follows that $H_3(G)$ is not finitely generated. So $G$ is not of type $FP_3$. Similarly, the result for $A*_C$ follows from the analogous sequence for HNN-extensions. □

In order to avail ourselves of Lemma 10 we need the following simple observation. (We continue with the notation established above.)

**Lemma 11.** *Suppose that $A$ has cohomological dimension 2 and that the Euler characteristic $\chi(Q)$ of $Q$ is zero. Then $H_2(C)$ is not finitely generated.*

*Proof.* Since $C$ is a finitely generated subgroup of a group of cohomological dimension 2, $H_2(C)$ is finitely generated if and only if $\chi(C) < \infty$. Suppose then that $\chi(C) < \infty$. Then

$$\chi(A) = \chi(Q)\chi(C).$$

But

$$\chi(A) = (|R| + 4|X|) - (|X| + 2) + 1 > 0,$$

which is a contradiction. □

In order to concoct some actual examples of finitely presented subgroups of automatic groups which are not automatic, notice that every two-generator one-relator group without torsion is of cohomological dimension at most two and has Euler characteristic zero. In particular, our constructions apply to the infinite cyclic group, to the free abelian group of rank two and to the Klein bottle group $\langle x, y \mid x^2 = y^2 \rangle$. Let us record a presentation of one such group.

**Example 1:** The following is a presentation for the double $G$ of a Rips group $A$ with quotient $\mathbb{Z}$ along its kernel $C$ which is the subgroup generated by $a$ and $b$. By Theorem 1 it is a subgroup of $F \times A$, where $F$ is a free group of rank two generated by $x$ and $\bar{x}$. $G$ is not of type $FP_3$.

$$\begin{aligned}
G = < x, \bar{x}, a, b | \ x^{-1}ax &= (ab)^{81}(ab^2)^{82}\cdots(ab)^{159}(ab^2)^{160}\\
xax^{-1} &= (ab)^{161}(ab^2)^{162}\cdots(ab)^{239}(ab^2)^{240}\\
x^{-1}bx &= (ab)^{241}(ab^2)^{242}\cdots(ab)^{319}(ab^2)^{320}\\
xbx^{-1} &= (ab)^{321}(ab^2)^{322}\cdots(ab)^{399}(ab^2)^{400}\\
\bar{x}^{-1}a\bar{x} &= (ab)^{81}(ab^2)^{82}\cdots(ab)^{159}(ab^2)^{160}\\
\bar{x}a\bar{x}^{-1} &= (ab)^{161}(ab^2)^{162}\cdots(ab)^{239}(ab^2)^{240}\\
\bar{x}^{-1}b\bar{x} &= (ab)^{241}(ab^2)^{242}\cdots(ab)^{319}(ab^2)^{320}\\
\bar{x}b\bar{x}^{-1} &= (ab)^{321}(ab^2)^{322}\cdots(ab)^{399}(ab^2)^{400} >
\end{aligned}$$

$\square$

We shall return to this example below, where we show that its isoperimetric function is strictly exponential.

**4. Examples of Stallings' type.** Here we concoct some examples of doubles in which the subgroup along which we double is not finitely generated but the resulting group is. In order to ensure that the groups obtained in this way are not asynchronously automatic, we must double twice.

To begin with we apply Theorem 1 in the case where $A$ is a free group of rank two with basis $x$ and $y$, and $C$ is the normal subgroup generated by all of the conjugates of $y$ by the powers of $x$. Let $G$ be the double of $A$ along $C$:

$$G = \langle x, y, \bar{x}, \bar{y} \mid x^{-i}yx^i = \bar{x}^{-i}\bar{y}\bar{x}^i \ (i \in \mathbb{Z})\rangle.$$

It follows from the Mayer-Vietoris sequence that $H_2(G)$ is not finitely generated. Since $d(A/C) = 1$, Theorem 7 provides us with a finitely presented group $E$ which can be presented in the form:

$$E = \langle x, y, \bar{x}, \bar{y}, s \mid s^{-1}xs = x, s^{-1}ys = x^{-1}yx, s^{-1}\bar{x}s = \bar{x},$$
$$s^{-1}\bar{y}s = \bar{x}^{-1}\bar{y}\bar{x}, x^{-1}yx = \bar{x}^{-1}\bar{y}\bar{x}\rangle.$$

Notice that $G$ is a normal subgroup of $E$. We now form $S$, the double of $E$ along $G$. This is the group of Stallings to which we alluded in the introduction. It follows from the Mayer-Vietoris sequence that $H_3(S)$ is not finitely generated. According to Theorem 7, $E$ embeds in $(F \times \mathbb{Z}) \times A$, where $F$ denotes a free groups of rank two, and by Theorem 1 $S$ embeds in $F \times E$. In other words, $S$ embeds in

$$F \times (F \times \mathbb{Z}) \times F.$$

This is our version of Mess' observation concerning Bieri's embedding of Stallings' group. Unlike our examples of Rips' type, Stallings' groups satisfies a polynomial isoperimetric inequality. Indeed, below we give an easy argument to show that it satisfies a cubic isoperimetric inequality; we do not know if this is optimal.

For the convenience of the reader, we provide a presentation of $S$:

$$S = \langle x, y, \bar{x}, \bar{y}, s, t \mid s^{-1}xs = x, s^{-1}ys = x^{-1}yx, s^{-1}\bar{x}s = \bar{x}, \ s^{-1}\bar{y}s = \bar{x}^{-1}\bar{y}\bar{x},$$
$$t^{-1}xt = x, \ t^{-1}yt = x^{-1}yx, \ t^{-1}\bar{x}t = \bar{x}, \ t^{-1}\bar{y}t = \bar{x}^{-1}\bar{y}\bar{x}, \ y = \bar{y}\rangle.$$

This may be expressed more concisely as:

$$S = \langle a, b, c, d, e \mid a^b = a^c = a^d = a^e, [c,d] = [d,b] = [e,c] = [e,b] = 1 \rangle.$$

For our next example, we again choose $A$ to be the free group on $x$ and $y$. However we now take $C$ to be the commutator subgroup $F'$ of $F$. $C$ is an infinitely generated subgroup of $F$. Indeed it is freely generated by the elements $[x,y]^{x^i y^j}$, with $i$ and $j$ ranging over all of the integers. As usual now, we form the double $G$ of $A$ along $C$, and then form the semidirect product $E$ of $G$ with the free group of rank two on $s$ and $t$. $E$ is a finitely presented group, with presentation

$$E = \langle x, y, \bar{x}, \bar{y}, s, t \mid s^{-1}xs = x, s^{-1}ys = x^{-1}yx, s^{-1}\bar{x}s = \bar{x}, s^{-1}\bar{y}s = \bar{x}^{-1}\bar{y}\bar{x},$$
$$t^{-1}xt = y^{-1}xy, t^{-1}yt = y, t^{-1}\bar{x}t = \bar{y}\bar{x}\bar{y}, t^{-1}\bar{y}t = \bar{y}, [x,y] = [\bar{x},\bar{y}] \rangle.$$

Notice that $E$ embeds in

$$H = (\mathbb{Z}^2 * \mathbb{Z}^2) \times F \times F.$$

We now repeat the process and form

$$K = E *_G \bar{E}.$$

Then $K$ turns out to be a subgroup of a direct product of a free group of rank four and the group $H$ described above. Notice that as usual, $H_3(K)$ is not finitely generated and so $K$ is not asynchronously automatic. It is not difficult to show that $E$ satisfies a polynomial isoperimetric inequality; we omit the details.

**5. Subgroups with exponential isoperimetric functions.** In this section we prove a result concerning the distortion of finitely generated normal subgroups of hyperbolic groups, and use this to concoct an abundance of subgroups of biautomatic groups whose isoperimetric functions are strictly exponential.

**Theorem 2** *Let $1 \to C \to A \to Q \to 1$ be a short exact sequence of finitely generated infinite groups. Suppose that $A$ is word hyperbolic and $Q$ is biautomatic. Then, the double $G$ of $A$ along $C$ satisfies, at best, an exponential isoperimetric inequality.*

It is important to note that Theorem 1 implies that the groups $G$ yielded by Theorem 2 are subgroups of biautomatic groups. Thus we obtain the first examples of subgroups of automatic groups for which the optimal isoperimetric inequality is exponential.

Before proceeding to the proof of the theorem we consider some explicit situations to which it applies.

**Example $1'$: Rips' type.** Applying Theorem 2 to Rips groups with biautomatic quotients yields many examples of subgroups of biautomatic groups that are not of type $FP_3$ and which have strictly exponential isoperimetric functions. In particular, the group considered in Example 1 above has this property. □

**Example 2: Bestvina-Feighn type.** In [4] Bestvina and Feighn studied semidirect products of the form $A = C \rtimes_\varphi \mathbb{Z}$, where $C$ is a finitely generated free group. They proved that such groups are hyperbolic if and only if they do not contain a free abelian subgroup of rank two. We assume that this is the case. Notice that in this situation the double of $A$ along $C$ is $G \cong C \rtimes F$, where $F$ is a free group of rank two on $s$ and $\bar{s}$ and both $s$ and $\bar{s}$ act on $C$ by the automorphism $\varphi$. According to Theorem 1, we may embed $G$ in the biautomatic group $F \times A$. According to Theorem 2, the isoperimetric function of $G$ is at best exponential; in fact $G$ is asynchronously automatic [6], so its isoperimetric function is strictly exponential.

An automorphism $\varphi$ of a finitely generated free group $C$ is said to have a periodic conjugacy class if there exist a non-zero integer $m$ and a pair of elements $v, w \in C$, $w \neq 1$, such that $\varphi^m(w) = v^{-1}wv$. It is not difficult to show that $C \rtimes_\varphi \mathbb{Z}$ has no free abelian subgroup of rank two if and only if $\varphi$ does not have any periodic conjugacy classes. Large classes of such automorphisms are known. This then provides us with a way of making explicit groups to which we can easily apply Theorems 1 and 2.

The following example is due to Gersten and Stallings [14]. Let $C$ be a free group of rank 3 with generators $a, b, c$ and let $\psi$ be the automorphism

$$a \mapsto c, \ b \mapsto ac, \ c \mapsto bc.$$

This is a so-called $PV$-automorphism, and as such has no periodic conjugacy classes. The double of $A = C \rtimes_\psi \mathbb{Z}$ along $C$ has presentation

$$\langle a, b, c, s, \bar{s} \mid s^{-1}as = c = \bar{s}^{-1}a\bar{s}, \ s^{-1}bs = ac = \bar{s}^{-1}b\bar{s}, \ s^{-1}cs = bc = \bar{s}^{-1}c\bar{s} \rangle.$$

□

**Example 3: Hyperbolic 3-manifolds.** There are many closed hyperbolic 3-manifolds which fibre over the circle. Indeed Thurston has conjectured that every closed hyperbolic 3-manifold has a finite sheeted cover that fibres in this way [19]. If $M^3$ admits such a fibration then its fundamental group $A = \pi_1 M^3$ fits into a short exact sequence $1 \to C \to A \to \mathbb{Z} \to 1$, where $C$ is the fundamental group of a closed surface of genus at least two. Applying Theorem 2 to this situation we deduce that $G$, the double of $A$ along the fibre group $C$, has a strictly exponential isoperimetric function. In this case $G$ is asynchronously automatic [6].

The construction of Theorem 1 embeds $G$ in the biautomatic group $\pi_1 M^3 \times F$, where $F$ is a free group of rank two. This group in turn embeds in $\Gamma = \pi_1 M^3 \times \pi_1 \Sigma$, where $\Sigma$ is any closed surface of genus at least two. Thus the group $G$, which has a strictly exponential isoperimetric function, embeds in $\Gamma$, which is the fundamental group of a closed 5-manifold of non-positive curvature. Notice that $\Gamma$ is a direct product of hyperbolic groups and hence is biautomatic. □

Combining this consequence of Theorem 2 with an earlier example we can easily deduce Theorem 3.

*Proof of Theorem 3.* Let $G$ and $\Gamma$ be the groups constructed in Example 3. Let $S$ be the example of Stallings type constructed in the previous section for which $H_3(S)$ is infinitely generated. Since $S$ embeds in a direct product of four free groups, it also embeds in a direct product of four surface groups $\pi_1 \Sigma$ of genus at least two. Hence $G \times S$ embeds in $B = \Gamma \times (\pi_1 \Sigma)^4$ which is the fundamental group of a manifold of non-positive curvature and is biautomatic. But $G \times S$ has a strictly exponential isoperimetric function and $H_3(G \times S)$ is not finitely generated. This completes the proof. □

Our proof of Theorem 2 relies upon the following lemma. We refer to [1] for basic facts about word hyperbolic groups. Recall that if $\Gamma$ is a finitely generated group with word metric $d$, then a cyclic subgroup $\langle x \rangle$ is said to be *undistorted* if there exists a constant $\lambda > 0$ such that $d(x^n, x^m) \geq \lambda |n - m|$ for all integers $n$ and $m$.

**Lemma 12.** *Let $1 \to C \to A \to Q \to 1$ be an exact sequence of finitely generated groups; let $d_C$ be the word metric corresponding to a choice of finite generating set for $C$. Suppose that $A$ is word hyperbolic and that $Q$ contains an undistorted cyclic subgroup $\langle x \rangle$. Let $a$ be an element in the preimage of $x$. Then, for every element $c \in C$ of infinite order there exist constants $\alpha > 1$ and $\beta > 0$ such that $d_C(1, a^{-n}ca^n) > \alpha^n - \beta$ for every $n \in \mathbb{N}$.*

*Proof.* Let $c, x$ and $a$ be as in the statement of the lemma. Let $\lambda$ be as in the sentence preceding the statement of the lemma. We choose a finite set of generators for $Q$, we then choose one element of $A$ from the preimage of each of these generators and append these elements to our chosen set of generators for $C$ to obtain a system of generators for $A$. With respect to the resulting word metrics, $A \to Q$ does not increase distances. Since $\langle x \rangle$ is undistorted in $Q$, it follows that the maps $\mathbb{N} \to A$ given by $n \mapsto a^n$ and $n \mapsto ca^n$ are quasigeodesics in $A$, indeed $d(1, a^n) \geq \lambda n$ for all $n$.

A pair of quasigeodesic rays in a word hyperbolic group either lie within a bounded distance of one another or else diverge exponentially. If $n \mapsto a^n$ and $n \mapsto ca^n$ were boundedly close then there would exist a $K > 0$ such that for each $n \geq 0$, there is a $m \geq 0$ such that $d_A(a^m, ca^n) = d_A(1, a^{-m}ca^n) < K$. There are only finitely many points $a^{-m}ca^n$ in any ball of finite radius about the identity in $A$, so such a bound would imply that there exist integers $m, n, p$ and $q$ with $m \neq p$, $n \neq q$ such that $a^{-m}ca^n = a^{-p}ca^q$, whence $c^{-1}a^{m-p}c = a^{n-q}$. The fact that the image of $a$ has infinite order in $Q$ implies that in this case $c$ would be in the centralizer of a power of $a$. But the centralizer of an element of a hyperbolic group is virtually cyclic [16]. Thus some power of $c$ would be equal to some power of $a$, contradicting the fact that $c$ has infinite order in $A$ and the image of $a$ has infinite order in $Q$ (the image of $c$ is the identity element). It follows that the rays $n \mapsto a^n$ and $n \mapsto ca^n$ must diverge exponentially.

We need to explain this last assertion more carefully: there exist geodesic rays $\gamma$ and $\gamma'$ beginning at the identity that are Hausdorff close to $n \mapsto a^n$ and $n \mapsto ca^n$ respectively. A simple exercise in the triangle inequality shows that there are positive constants $\mu$ and $\nu$ such that for every positive integer $n > 1/\lambda$ there exist integers $m_1$ and $m_2$ with $n \geq m_2 \geq \lambda n$ and $\lambda n > m_1 \geq \nu \lambda n$ and $d(\gamma(m_2), a^n) \leq \mu$ and $d(\gamma'(m_2), ca^n) \leq \mu$. Because $A$ is hyperbolic there exists a constant $\zeta > 1$ such that if $p$ is path in the Cayley graph of $A$ that joins $\gamma(m)$ to $\gamma'(m)$ and if this path lies in the closure of the complement of the ball of radius $m$ about the identity, then $p$ has length at least $\zeta^{m-M}$, where $M$ is a constant depending on $\gamma$ and $\gamma'$ (see for example [1], 2.19).

In order to complete the proof of the lemma, we must exhibit an exponential lower bound on the length of any path in the Cayley graph of $A$ that joins $a^n$ to $ca^n$ and is labelled entirely by generators of $C$. Let $w_n$ be the word labelling such a path. Consider the path joining $\gamma(m_1)$ to $\gamma'(m_1)$ that first traverses the segment $[\gamma'(m_1), \gamma'(m_2)]$ of $\gamma'$, then follows a geodesic from $\gamma'(m_2)$ to $ca^n$, then follows the path labelled $w_n$ from $ca^n$ to $a^n$, then follows a geodesic from $a^n$ to $\gamma(m_2)$, then traverses the segment $[\gamma(m_2), \gamma(m_1)]$ in the image of $\gamma$. There is an obvious linear bound, $\xi n$ say, on the sum of all portions of this path except that labelled $w_n$. Thus, if we can show that this portion of the path does not enter the ball of radius $m_1$ about the identity, then we will obtain a lower bound of $\zeta^{m_1-M} - \xi n$ on the length of $w_n$. Since $m_1 \geq \nu \lambda n$, this will complete the proof of the lemma.

Let $b$ be a vertex (element of $A$) along the subpath labelled $w_n$, and let $u$ be the word labelling some geodesic from the identity to $b$. Note that $b = a^n c_b$ for some $c_b \in C$, and hence the image of $b$ in $Q$ is $x^n$. It follows that the number of edges of $u$ that are labelled by generators which project non-trivially to $Q$ is at least $d_Q(1, x^n) \geq \lambda n$. But $m_1 < \lambda n$. Since $u$ was geodesic, we deduce that $d_A(1, b) > m_1$, as required. □

*Proof of Theorem 2.* Let $C, A, Q$ and $G$ be as in the statement of the theorem. Every infinite subgroup of a hyperbolic group contains an element of infinite order; every infinite cyclic subgroup of a biautomatic group is undistorted [13]; and every infinite automatic group contains an element of infinite order [15]. Therefore we may choose elements $c \in C$ and $x \in Q$ (with lift $a \in A$) to which Lemma 12 applies.

We choose finite generating sets $x_1, \ldots, x_n$ for $Q$ and $c_1, \ldots, c_m$ for $C$, which include $x$ and $c$ respectively. Then we choose one element $a_i \in A$ in the preimage of each $x_i$, and work with the generating set $a_1, \ldots, a_n, \bar{a}_1, \ldots, \bar{a}_n, c_1, \ldots, c_m$ for $G$. We fix a finite presentation $\mathcal{P}$ for $G$ with this set of generators.

Consider the word $a^{-n} c a^n \bar{a}^{-n} c^{-1} \bar{a}^n$, which represents the identity in $G$. Let $\Delta_n$ be a van Kampen diagram over $\mathcal{P}$ for this word. Let $p_0$ denote the basepoint of the diagram and let $p_1$ be the vertex of $\partial \Delta$ corresponding to the subword $a^{-n} c a^n$. We wish to establish an exponential lower bound on the number of 2-cells in $\Delta$. In the light of the preceding lemma, it suffices to show that there is an injective path in the 1-skeleton of $\Delta$ that connects $p_0$ to $p_1$ and consists entirely of edges labelled by generators of the form $c_i$. (Lemma 12 gives an exponential lower bound on the length of such a path, and since it is injective and there is a uniform bound on the number of edges in the boundary of any 2-cell, we obtain an exponential lower bound on the number of distinct 2-cells in $\Delta_n$ that this path touches.)

Associated to the sequence

$$1 \to C \to G \to Q * \bar{Q} \to 1.$$

we have a continuous cellular map from the 1-skeleton $\Delta_n^{(1)}$ of $\Delta_n$ to the Cayley graph of $Q * \bar{Q}$. This map sends $p_0$ to the identity vertex, sends edges labelled $a_i$ homeomorphically onto edges labelled $x_i$, sends edges labelled $\bar{a}_i$ homeomorphically onto edges labelled $\bar{x}_i$, and restricts to a constant map on each edge labelled $c_j$. The inverse image of the identity vertex in the Cayley graph of $Q * \bar{Q}$ contains the subgraph of $\Delta_n^{(1)}$ spanned by those vertices which are connected to $p_0$ by a path whose edges are labelled by the $c_1, \ldots, c_m$. Let us denote this subgraph $\Gamma$. Notice that the generators labelling the edges of $\Gamma$ are of the form $c_j$.

By removing any vertex from the Cayley graph of $Q * \bar{Q}$ one disconnects it. Also, the two edges of $\partial \Delta_n$ that are incident at $p_0$ are mapped to distinct components of the Cayley graph of $Q * \bar{Q}$ minus the identity. Therefore, these two edges must lie in different connected components of the complement of $\Gamma$ in $\Delta_n^{(1)}$. But since $\Gamma$ can only meet $\partial \Delta_n$ at $p_0$ and $p_1$, this means that there must be a path in $\Gamma$ connecting $p_0$ to $p_1$. Any choice of a shortest such path yields the injective path required to complete the proof. □

**6. Closing remarks.** We shall now explain why Stallings' group $S$ satisfies a cubic isoperimetric inequality. We work with the second (more concise) presentation given above. Observe that $S$ can be viewed as the HNN-extension of its subgroup $H = gp\{b,c\} \times gp\{d,e\}$, which is a direct product of free groups, obtained by adding the stable letter $a$ which commutes with the subgroup consisting of all words with exponent sum zero.

Consider a word $w = w(a,b,c,d,e)$ in the generators $a,b,c,d,e$ and their inverses that represents the identity in $S$. Notice that all generators appear with exponent sum zero. Let $n$ denote the length of $w$. There is a constant $k$ such that by applying at most $kn^2$ relators of the form $[x,y] = 1$ we can transform $w$ so as to ensure that it has a subword of the form $a^\epsilon u_1(b,c) v_1(d,e) a^{-\epsilon}$ where the exponent sums in the various generators satisfy the equation $\sigma_b(u_1) + \sigma_c(u_1) + \sigma_d(v_1) + \sigma_e(v_1) = 0$, with $\epsilon = \pm 1$. For simplicity we assume that $\epsilon = 1$. By applying three relators we can perform the transformations

$$ab = b(b^{-1}ab) = ba^b = ba^e$$

and

$$a^e b = e^{-1}abe = e^{-1}ba^b e = ba^{e^2}.$$

In this way we can push the letter $a$ past the subword $u_1(b,c)$ by applying $K |u_1|^2$ relators, where $K$ is a universal constant. In the same way, one can pass $a^{-1}$ to the left through $v_1$, at quadratic cost. After having done so, the above word is transformed into

$$u_1 a^{e^p} a^{-b^p} v_1.$$

A further $|2p|^2$ applications of the relator $[e,b] = 1$ transform the subword $a^{e^p} a^{-b^p}$ to

$$e^{-p} a(eb^{-1})^p a^{-1} b^p.$$

Then $|p|$ uses of the relator $a^e = a^b$ transforms this to $e^{-p}(eb^{-1})^p b^p$ eliminating $a$. Then a further $|2p|^2$ transforms this to $e^{-p} e^p b^{-p} b^p$ which is freely equal to 1. So these steps transform the original subword to $u_1 v_1$. Combining these estimates we see that we have succeeded in reducing the number of occurrences of $a$ in $w$ at quadratic cost. Therefore, there exists a constant $\lambda > 0$ such that by applying at most $\lambda n^3$ relators we can remove all occurrences of $a$. But then, the remaining word, which is no longer than the original word $w$, represents the identity in $gp\{b,c\} \times gp\{d,e\}$, which is a product of free groups, and hence satisfies a quadratic isoperimetric inequality.

We end with an example that hints at the many possible variations on our earlier examples.

**Example 4.** Let $C$ denote the free group of rank three on $x_1, x_2, x_3$ and consider the automorphism $\varphi$ of $C$ given by

$$\varphi(x_1) = x_1, \ \varphi(x_2) = x_2 x_1, \ \varphi(x_3) = x_3 x_2.$$

The group $A = C \rtimes_\varphi \mathbb{Z}$ is biautomatic and $G$, its double along $C$, embeds in the biautomatic group $A \times F$, where $F$ is free of rank two. For every $c \in C$, the function $n \mapsto d(1, \varphi^n(c))$ is dominated by a quadratic function, and for some choices of $c$ (for example $c = x_3$) this function is $O(n^2)$. By using the techniques introduced in [7], one can deduce from this growth estimate that the isoperimetric function for $G$ is polynomial of degree 4. In a similar manner, one can construct polynomials of other degrees. □

Department of Mathematics, City College of New York, Convent Avenue at 138th Street, New York, NY 10031, USA

Mathematical Institute, 24–29 St. Giles, Oxford OX1 3LB, U.K.

Department of Mathematics, University of Melbourne, Parkville 3052, Australia

Centre de Mathématiques et d'Informatique, Université de Provence, F–13453 Marseille cedex 13, France